\documentclass{amsart}
\usepackage{amsmath,amssymb,amsfonts,amsthm}

\newtheorem{thm}{Theorem}[section]
\newtheorem{cor}[thm]{Corollary}
\newtheorem{lem}[thm]{Lemma}

\theoremstyle{definition}

\theoremstyle{remark}

\numberwithin{equation}{section}

\begin{document}

\title
 { On the Hausdorff dimension faithfulness connected with $Q_{\infty}$-expansion}

\author{Jia Liu  \& Zhenliang Zhang}


\address{Institute of Statistics and Applied Mathematics, Anhui University of Finance and Economics 233030,
Bengbu, P.R.China}

\email{liujia860319@163.com}
\address{
School of Mathematical Sciences, Henan Institute of Science and Technology, Xinxiang, 453003, P. R. China;
}
\email{zhliang\_zhang@163.com}

\date{}
\begin{center}

\end{center}


\begin{abstract}
In this paper, we show that,
the family of all possible union of finite consecutive cylinders
 of the same rank of $Q_{\infty}$-expansion is faithful for the Hausdorff dimension calculation.
Applying this result,
we give the necessary and sufficient condition for the family of all cylinders
 of $Q_{\infty}$-expansion to be faithful for Hausdorff dimension calculation on the unit interval,
this answers the open problem mentioned in \cite{AKNT2}.

\medskip

\smallskip


\end{abstract}
\maketitle

\section{Introduction}

The notion Hausdorff dimension is well-know now and plays an important role in fractal geometry.
Hausdorff dimension has the advantage of being defined any set,
 as it is based on measures,
 which are relatively easy to manipulate.
 However a major disadvantage is that in may cases
 it is a rather non-trivial problem to give the exact Hausdorff dimension.
 Recent years,
 there has been a large interest in determining the Hausdorff dimension
 by restricting the family of admissible coverings.
 In \cite{AKNT1,AKNT2},
 they introduce the notion of faithfulness and non-faithfulness of families for Hausdorff dimension calculation.
 The faithfulness of family leads us to calculate the Hausdorff dimension
 only considering narrow family of admissible coverings.

 To be more precise,
 let us shortly recall some definitions.
 Let $\Phi$ be a fine family of coverings on $[0,1]$,i.e.,
 a family of subsets of $[0,1]$ such that for any $\epsilon>0$,
 there exists an at most countable $\epsilon$-covering $\{E_j\}$ of $[0,1]$ with $E_j\in\Phi$.
 The $\alpha$-dimensional Hausdorff measure of a set $E\subset[0,1]$
 with respect to a fine family of coverings $\Phi$ is defined by
 $$
 H^{\alpha}(E,\Phi)=\liminf_{\epsilon\to 0}\left\{\sum_{j}|E_j|^{\alpha}|\text{ }\{E_j\}
 \text{ is an } \epsilon-\text{covering of E}, E_j\in\Phi\right\}
 $$
 and the nonnegative number
 $${\dim}_H(E,\Phi)=\inf\left\{\alpha|H^{\alpha}(E,\Phi)=0 \right\}.$$
is called the Hausdorff dimension of the set $E\subset[0,1]$ with respect to the family $\Phi$.
 If we take $\Phi$ as the family of all subsets of $[0,1]$,
 we denote ${\dim}_H(E,\Phi)$ by ${\dim}_H(E)$,
 which is equal to the classical Hausdorff dimension of set $E\subset[0,1]$.
For more properties of Hausdorff dimension, one is referred to \cite{KFM,Rog}.

 A fine covering family $\Phi$ is said to be a faithful family of coverings (non-faithful family of coverings)
 for Hausdorff dimension calculation if
 $${\dim}_H(E,\Phi)={\dim}_H(E), \text{ for any } E\in [0,1]$$
 $$(\text{respect to } \exists E\subset [0,1]:{\dim}_H(E,\Phi)\neq{\dim}_H(E) ).$$

The first result concerning the problem of faithfulness coverings is due to A. Besicovitch \cite{Bec},
who proved that the the family of cylinders of a binary expansion.
This result was extended to the family of s-adic cylinders by P. Billingsley\cite{Bil}.
S. Albeverio, G. Ivanenko, M. Lebid, G. Torbin et al. have done a series of work in this direction,
e.g.  \cite{AKNT1,AKNT2,APT,AT,PT}.
Especially, in \cite{AKNT2},
S. Albeverio, Y. Kondratiev, R. Nikiforov and G. Torbin \cite{AKNT2} gave a general sufficient conditions
for the family of all cylinders of ${Q_{\infty}}$-expansion to be faithful.

In this paper,
we focus on the family of the finite union of cylinders of ${Q_{\infty}}$-expansion
and prove that the family of all possible union of finite consecutive cylinders
of the same rank of $Q_{\infty}$-expansion is faithful for the Hausdorff dimension calculation.
Applying this result,
we give the necessary and sufficient condition for the family
of all cylinders of $Q_{\infty}$-expansion to be faithful family for Hausdorff dimension calculation
on the unit interval, this answers the open problem mentioned in \cite{AKNT2}.

\section{Statement of main results}

First of all,
we briefly recall the definition and some properties related to $Q_{\infty}$-expansion.

Let $Q=(q_0,q_1,\cdots,q_n,\cdots,)$ be a $Q_{\infty}$-vector,
i.e., $\{q_i\}_{i\geq 0}$ are positive reals and satisfy $\sum_{i\geq 0}q_i=1$.
Now we give the $Q_{\infty}$-expansion with respect to the $Q_{\infty}$-vector in the following way.

Step 1.
We decompose $[0,1)$ (from the left to the right) into the union of semi-intervals
$\Delta_{\alpha_1},\alpha_1\in\mathbb{N}$ (without common points)
with length $|\Delta_{\alpha_1}|=q_{\alpha_1}$,
$$
[0,1)=\bigcup_{{\alpha_1}=0}^{\infty}\Delta_{\alpha_1}.
$$
Each interval $\Delta_{\alpha_1}$ is called a 1-rank cylinder.

Step $n\geq 2$.
We decompose (from the left to the right) each $(n-1)$-rank cylinder (semi-intervals)
$\Delta_{\alpha_1\cdots\alpha_{n-1}}$ into the union of $n$-rank cylinders
$\Delta_{\alpha_1\cdots\alpha_{n-1}\alpha_n}$ (without common points),
i.e.,
$$
\Delta_{\alpha_1\cdots\alpha_{n-1}}=\bigcup_{\alpha_{n}=0}^{\infty}\Delta_{\alpha_1\cdots\alpha_n}.
$$
Every $n$-rank cylinder has the length
$$
|\Delta_{\alpha_1\cdots\alpha_n}|=\prod_{i=1}^{n}q_{\alpha_i}.
$$

It is clear that any sequence of indices $\{\alpha_n\}$ generates
the corresponding sequence of embedded cylinders
$$
\Delta_{\alpha_1}\supset\Delta_{\alpha_1\alpha_2}\supset\cdots
\supset\Delta_{\alpha_1\alpha_2\cdots\alpha_n}\supset\cdots
$$
and there exists an unique points $x\in[0,1)$ belongs to all of them.

Conversely,
for any point $x\in[0,1)$ there is an unique sequence of embedded cylinders
$$
\Delta_{\alpha_1}\supset\Delta_{\alpha_1\alpha_2}\supset\cdots
\supset\Delta_{\alpha_1\alpha_2\cdots\alpha_n}\supset\cdots
$$
containing the point $x$,
i.e.,
$$
x=\bigcap_{n=1}^{\infty}\Delta_{\alpha_1\cdots\alpha_n}=
\bigcap_{n=1}^{\infty}\Delta_{\alpha_1(x)\cdots\alpha_n(x)}:=
\Delta_{\alpha_1(x)\cdots\alpha_n(x)}\cdots.
$$
The expression is called the $Q_{\infty}$-expansion of $x$.
Real numbers which are the endpoints of $n$-th rank cylinders are said to be $Q_{\infty}$-rational,
and their $Q_{\infty}$-expansion contains only finitely many non-zero digits.
\smallskip

Let $\Phi$ be the family of all possible cylinders of $Q_{\infty}$-partition of semi-interval $[0,1)$,
i.e.,
$$\Phi=\Phi(Q)=\{\Delta_{\alpha_1\cdots\alpha_n}:\alpha_i\in\mathbb{N},i=1,2,\cdots,n; n\in\mathbb{N}\}.$$

 In recent paper\cite{AKNT2}, S.Albeverio,Y.Kondratiev,R.Nikiforov and G.Torbin gave a general sufficient condition for the family $\Phi$ to be faithful.  However, based on the method which was invented by Yuval Peres, they proved the following result shows that the family $\Phi$ is not necessarily faithful.
 \begin{thm}\cite{AKNT2}
 If there exist constants $m_0>1$, $A>0$ and $B>0$ such that $$\frac{A}{i^{m_0}}\leq q_i\leq\frac{B}{i^{m_0}},\forall i\in\mathbb{N},$$
 then the family $\Phi$ is non-faithful.

 \end{thm}

In this paper,
let $\mathcal{A}_n$ be the family of all possible union of finite consecutive $n$-th rank cylinders,
i.e.,
$$\mathcal{A}_n=\left\{\bigcup_{i=m}^{m+k}\Delta_{\alpha_1\cdots\alpha_{n-1}i}:m,k,\alpha_j\in\mathbb{N},1
\leq j\leq n-1\right\}$$
and $$\mathcal{A}=\mathcal{A}(Q)=\bigcup_{n\geq 1}\mathcal{A}_n.$$
We mainly consider the covering family $\mathcal{A}$ of all possible union of finite cylinders,
and show that $\mathcal{A}$ is faithful for Hausdorff dimension calculation
without any additional condition with respect to $Q_{\infty}$-vector.
 \begin{thm} \label{Mthm1}
 Let $Q=(q_0,q_1,\cdots,q_n,\cdots,)$ be a $Q_{\infty}$-vector.
 Then the family $\mathcal{A}=\mathcal{A}(Q)$ is faithful for the Hausdorff dimension calculation on the unit interval.
 \end{thm}


In fact,
$Q_{\infty}$-expansion may be generated by some infinite linear IFS $\{F_0,F_1,\break\cdots,F_n,\cdots\}$
such that $F_n$ is a similar transformation with ratio $q_n$ and $\{sup F_n([0,1])\}$ is strictly monotone.
If we take $q_i=\frac{1}{(i+1)(i+2)}$,
then we get the following corollary for the classical L\"{u}roth expansion\cite{DK}.

\begin{cor}
The family of all possible union of finite consecutive cylinders in the same rank
generated by $L\ddot{u}roth$ expansion is faithful.
\end{cor}

Applying Theorem \ref{Mthm1}
we give the necessary and sufficient condition for the family $\Phi$ to be faithful
for the Hausdorff dimension calculation,
which answers the open problem mentioned in \cite{AKNT2}.
 \begin{thm} \label{Mthm2}
 Let $Q=(q_0,q_1,\cdots,q_n,\cdots,)$ be a $Q_{\infty}$-vector. Then the family $\Phi=\Phi(Q)$ is faithful for the Hausdorff dimension calculation on the unit interval if and only if for any $\alpha\in(0,1)$, $\delta\in(0,\alpha)$ there exists a positive integer $N=N(\alpha,\delta)$ such that for any $n,M>N$,
 \begin{equation}\label{equ1}
 \left(\sum_{i=n}^{n+M}q_i\right)^{\alpha-\delta}\geq\sum_{i=n}^{n+M}q_{i}^{\alpha}.
 \end{equation}
 \end{thm}

\section{Proofs of main results}
Before the proof of the main results, we give the following useful lemma, which can be proved easily.
\begin{lem}\label{lem1}
Let $\{a_n\}_{n\geq 0}$ be a sequence of positive reals,
if $\sum_{n\geq 0}a_n<\infty$,
then for any $\alpha\in(0,1]$
there exists a sequence $\{n_k\}_{k\geq 1}$ of integers with
$$
\left(\sum_{i=0}^{n_1}a_i\right)^{\alpha}\geq\sum_{i=1}^{\infty}\left(\sum_{j=n_{i}+1}^{n_{i+1}}a_j\right)^{\alpha}.
$$
\end{lem}
\noindent \textbf{Proof of the Theorem \ref{Mthm1}.}
We only need to show that $${\dim}_{H}(E)={\dim}_{H}(E,\mathcal{A})$$ for any set $E\subset[0,1]$.

Since the set of all $Q_{\infty}$-rational is dense in $[0,1]$,
for the calculation of the Hausdorff dimension of a set  $E\subset[0,1]$,
we may only consider the coverings of intervals with endpoints are $Q_{\infty}$-rational.

For a given set $E$,
$\alpha\in(0,1)$, $\delta\in(0,\alpha)$ and $\epsilon>0$,
let $\{E_j\}_{j\in\mathbb{N}}$ be an $\epsilon$-covering of $E$ with $E_j=[a_j,b_j)$
and $a_j,b_j$ are $Q_{\infty}$-rational.
For each interval $E_j, j\in\mathbb{N}$,
there exists a cylinder $\Delta_{\alpha_1\alpha_2\cdots\alpha_{n_j}}\in\Phi$ such that:

1) $E_j\subset\Delta_{\alpha_1\alpha_2\cdots\alpha_{n_j}}$;

2) any interval of rank $n_j+1$ does not contain $E_j$.

\noindent Without loss of generality,
we assume that $a_j$ have the following $Q_{\infty}$-expansion:
$$a_j=\Delta_{\alpha_1\cdots\alpha_{n_j}\cdots\alpha_{n_j+l_j}00\cdots}$$
for some integer $l_j$,
where $\alpha_k=\alpha_k(a_j)$ is the $k$-th digit of $a_j$.
Then the point $c_j=\Delta_{\alpha_1\cdots\alpha_{n_j}(\alpha_{n_j+1}+1)00\cdots}$ belongs to $[a_j,b_j)$.

For the simplicity of the notation,
take $\beta_k=\alpha_{n_j+k}(\alpha_j)$.
To cover $E_j$ by union of finite cylinders we consider the coverings of $[a_j,c_j)$ and $[c_j,b_j)$ separately.

\smallskip
\begin{itemize}
  \item The coverings of $[c_j,b_j)$.
\end{itemize}

If $b_j$ is the endpoint of some $(n_j+1)$-th rank cylinder, i.e.,
$b_j$ have the following $Q_{\infty}$-expansion:
$$b_j=\Delta_{\alpha_1\cdots\alpha_{n_j}(\beta_1+s_0)00\cdots}$$ for some integer $s_0\geq 2$.
Then $[c_j,b_j)$ is the union of finite cylinders of rank $n_j+1$,
i.e.,
\begin{equation}
[c_j,b_j)=\bigcup_{i=\beta_1+1}^{\beta_1+s_0-1}\Delta_{\alpha_1\cdots\alpha_{n_j-1}\alpha_{n_j}i}\in\mathcal{A}.
\end{equation}

If $b_j$ is not the endpoint of any $(n_j+1)$-th rank cylinder.
Then there exists an integer $s_1\in\mathbb{N}$ such that
$$b_j\in\text{Int}(\Delta_{\alpha_1\cdots\alpha_{n_j}(\beta_1+s_1+1)}),$$
where Int(A) denotes the interior of set $A$.
Let $s_2$ be the integer with
$$b_j\in\Delta_{\alpha_1\cdots\alpha_{n_j}(\beta_1+s_1+1)\underbrace{0\cdots 0}_{s_2}}$$
and
$$\Delta_{\alpha_1\cdots\alpha_{n_j}(\beta_1+s_1+1)\underbrace{0\cdots 0}_{s_2+1}}\subset[c_j,b_j).$$
Then we have $$|\Delta_{\alpha_1\cdots\alpha_{n_j}(\beta_1+s_1+1)\underbrace{0\cdots 0}_{s_2}}|=\frac{1}{q_0}\Delta_{\alpha_1\cdots\alpha_{n_j}(\beta_1+s_1+1)\underbrace{0\cdots 0}_{s_2+1}}\leq\frac{1}{q_0}|E_j|.$$

Take  $$J_0=\bigcup_{i=1}^{s_1}\Delta_{\alpha_1\cdots\alpha_{n_j}(\beta_1+i)} \text{ if } s_1\geq 1 \ (J_0=\O \text{ if }s_1=0)$$ and
$$J_1=\Delta_{\alpha_1\cdots\alpha_{n_j}(\beta_1+s_1+1)\underbrace{0\cdots 0}_{s_2}}.$$
So $[c_j,b_j)$ can be covered by $\{J_0,J_1\}\subset\mathcal{A}$
such that the corresponding $\alpha$-volume does not exceed
\begin{equation}
(1+\frac{1}{q_0^{\alpha}})|E_j|^{\alpha}.
\end{equation}
\smallskip

\begin{itemize}
  \item  The coverings of $[a_j,c_j)$.
\end{itemize}

We can divide the interval $[a_j,c_j)$ into the union of pairwise disjiont cylinders of the following ranks.
\[
\begin{split}
&\text{Rank } (n_j+2) \text{: }\ \  \ \Delta_{\alpha_1\cdots\alpha_{n_j}\beta_1i}, \  \  i\geq\beta_2+1,\\
&\text{Rank } (n_j+3) \text{: }\  \ \ \Delta_{\alpha_1\cdots\alpha_{n_j}\beta_1\beta_2i},\  \  i\geq\beta_3+1,\\
&\text{}\cdots\\
&\text{Rank } (n_j+l_j-1) \text{: }\  \ \ \Delta_{\alpha_1\cdots\alpha_{n_j}\beta_1\cdots\beta_{l_j-2}i},\ \  i
     \geq\beta_{l_j-1}+1,\\
&\text{Rank } (n_j+l_j)\text{: }\  \ \ \Delta_{\alpha_1\cdots\alpha_{n_j}\beta_1\cdots\beta_{l_j-1}i},\ \ i\geq\beta_{l_j}.
\end{split}
\]
For the cylinders of rank $n_j+k, 2\leq k\leq l_j-1,$
by Lemma \ref{lem1},
there exists a sequence  $\{t_k(0)=0,t_k(i)\}_{i\geq 1}$ of integers such that
$$
\left|\bigcup_{i=\beta_k+1}^{\beta_k+t_k(1)}\Delta_{\alpha_1\cdots\alpha_{n_j}
\beta_1\cdots\beta_{k-1}i}\right|^{\alpha}
\geq
\sum_{m=1}^{\infty}\left|\bigcup_{i=\beta_k+t_k(m)+1}^{\beta_k+t_k(m+1)}
\Delta_{\alpha_1\cdots\alpha_{n_j}\beta_1\cdots\beta_{k-1}i}\right|^{\alpha}
$$
and
$$
\left|\bigcup_{i=\beta_{l_j}}^{\beta_{l_j}+t_{l_j}(1)}
\Delta_{\alpha_1\cdots\alpha_{n_j}\beta_1\cdots\beta_{l_j-1}i}\right|^{\alpha}
\geq
\sum_{m=1}^{\infty}\left|\bigcup_{i=\beta_{l_j}+t_{l_j}(m)+1}^{\beta_{l_j}+t_{l_j}(m+1)}
\Delta_{\alpha_1\cdots\alpha_{n_j}\beta_1\cdots\beta_{l_j-1}i}\right|^{\alpha}.
$$
Take
$$
J_k(m)=\bigcup_{i=\beta_k+t_k(m)+1}^{\beta_k+t_k(m+1)}
\Delta_{\alpha_1\cdots\alpha_{n_j}\beta_1\cdots\beta_{k-1}i},
\ 2\leq k\leq l_j-1, m\in\mathbb{N}
$$
and
\[
\begin{split}
J_{l_j}(m)&=\bigcup_{i=\beta_{l_j}+t_{l_j}(m)+1}^{\beta_{l_j}+t_{l_j}(m+1)}
\Delta_{\alpha_1\cdots\alpha_{n_j}\beta_1\cdots\beta_{l_j-1}i},\ m\geq 1,\\
J_{l_j}(0)&=\bigcup_{i=\beta_{l_j}}^{\beta_{l_j}+t_{l_j}(1)}
\Delta_{\alpha_1\cdots\alpha_{n_j}\beta_1\cdots\beta_{l_j-1}i}.
\end{split}
\]
Then
 $[a_j,b_j)$ can be covered by family $\{J_k(m): 2\leq k\leq l_j, m\in\mathbb{N}\}\subset\mathcal{A}$
 and corresponding $\alpha$-volume $V_j(\alpha)$,
\[
\begin{split}
V_j(\alpha)
&=\sum_{k=2}^{l_j}\sum_{m=0}^{\infty}|J_{k}(m)|^{\alpha}\leq 2\left[\sum_{k=2}^{l_j}|J_{k}(0)|^{\alpha}\right]\\
&= 2\left[\sum_{k=2}^{l_j-1}\left|\bigcup_{i=\beta_k+1}^{\beta_k+t_k(1)}
\Delta_{\alpha_1\cdots\alpha_{n_j}\beta_1\cdots\beta_{k-1}i}\right|^{\alpha}
+\left|\bigcup_{i=\beta_{l_j}}^{\beta_{l_j}+t_{l_j}(1)}
\Delta_{\alpha_1\cdots\alpha_{n_j}\beta_1\cdots\beta_{l_j-1}i}\right|^{\alpha}\right].
\end{split}
\]
Let $\mathbf{q}=\max_i(q_i)$ and $d_j\in\mathbb{N}$ such that $2^{d_j-1}<l_j\leq 2^{d_j}$,
then for any $2\leq k\leq l_j-1$,
\[
\begin{split}
\frac{|J_k(0)|^{\alpha}}{|E_j|^{\alpha-\delta}}
&=\frac{\left|\bigcup_{i=\beta_k+1}^{\beta_k+t_k(1)}
\Delta_{\alpha_1\cdots\alpha_{n_j}\beta_1\cdots\beta_{k-1}i}\right|^{\alpha}}{|E_j|^{\alpha-\delta}}\\
&\leq\frac{\left|\bigcup_{i=\beta_k+1}^{\beta_k+t_k(1)}
\Delta_{\alpha_1\cdots\alpha_{n_j}\beta_1\cdots\beta_{k-1}i}\right|^{\alpha-\delta}
\left|\Delta_{\alpha_1\cdots\alpha_{n_j}\beta_1\cdots\beta_{k-1}}\right|^{\delta}}{|E_j|^{\alpha-\delta}}\\
&\leq|\Delta_{\alpha_1\cdots\alpha_{n_j}\beta_1\cdots\beta_{k-1}}|^{\delta}\leq\mathbf{q}^{(k-1)\delta}
\end{split}
\]
and similarly, $$\frac{|J_{l_j}(0)|^{\alpha}}{|E_j|^{\alpha-\delta}}\leq\mathbf{q}^{(l_{j}-1)\delta},$$
thus we have the following estimation of $V_j(\alpha)$.
\[
\begin{split}
V_j(\alpha)
&\leq 2\left[\sum_{i=1}^{d_j-1}\sum_{2^{i-1}<k\leq 2^i}|J_k(0)|^{\alpha}+
\sum_{2^{d_j-1}<k\leq l_j}|J_k(0)|^{\alpha}\right]\\
&\leq\frac{2}{\mathbf{q}^{\delta}}|E_j|^{\alpha-\delta}
    \sum_{k=0}^{d_j-1}2^k(\mathbf{q}^{\delta})^{2^k}
    \leq \frac{2}{\mathbf{q}^{\delta}}|E_j|^{\alpha-\delta}
    \sum_{s=1}^{\infty}s(\mathbf{q}^{\delta})^{\frac{s}{2}}(\mathbf{q}^{\delta})^{\frac{s}{2}}\\
&\leq\frac{2}{\mathbf{q}^{\delta}}|E_j|^{\alpha-\delta}W(\delta)
\sum_{s=1}^{\infty}(\mathbf{q}^{\delta})^{\frac{s}{2}}=
  \frac{2W(\delta)}{(1-\mathbf{q}^{\frac{\delta}{2}})\mathbf{q}^{\frac{\delta}{2}}}|E_j|^{\alpha-\delta}
\end{split}
\]
where $W(\delta)$ is a constant such that $s(\mathbf{q}^{\delta})^{\frac{s}{2}}\leq W(\delta)$ for any $s\in \mathbb{N}$.

So,
for a given $E_j=[a_j,b_j)$
there exists a countable subfamily of $\mathcal{A}$ that cover $E_j$
and the corresponding $\alpha$-volume does not exceed
$$
K(\alpha,\delta)|E_j|^{\alpha-\delta}
$$
where
$K(\alpha,\delta)=1+\frac{1}{q_{0}^{\alpha}}
+\frac{2W(\delta)}{(1-\mathbf{q}^{\frac{\delta}{2}})\mathbf{q}^{\frac{\delta}{2}}}$.

Therefore, for any $E\subset[0,1]$, $\alpha\in(0,1)$, $\delta\in(0,\alpha)$ we have $$H^{\alpha}(E)\leq H^{\alpha}(E,\mathcal{A})\leq K(\alpha,\delta)H^{\alpha-\delta}(E),$$ this gives ${\dim}_{H}(E,\mathcal{A})\leq{\dim}_{H}(E)+\delta$ for any $\delta\in(0,\alpha)$, thus we show that $${\dim}_{H}(E)={\dim}_{H}(E,\mathcal{A}),$$
which completes the proof.

\bigskip

\noindent \textbf{Proof of the Theorem \ref{Mthm2}.}
\smallskip

\textbf{Sufficiency.}
Suppose inequality (\ref{equ1}) hold.
Due to Theorem \ref{Mthm1},
it is enough to consider the coverings coming from $\mathcal{A}$ for the calculation
of the Haudorff dimension of a set $E\subset[0,1]$.
Let $\epsilon>0$, $\{E_j\}\subset\mathcal{A}$ be a $\epsilon$-covering of $E$,
without loss of generality,
we assume that $E_j$ is an union of finite consecutive cylinders of rank $n_j$, i.e.,
$$
E_j=\bigcup_{i=m}^{m+k}\Delta_{\alpha_1\cdots\alpha_{n_j-1}i},
$$
for some $k,m,\alpha_1,\cdots,\alpha_{n_j-1}\in\mathbb{N}$.

If $k\leq N$,
then $$|E_j|^{\alpha}\geq\frac{1}{N}\sum_{i=m}^{m+k}|\Delta_{\alpha_1\cdots\alpha_{n_j-1}i}|^{\alpha}.$$

If $k>N$,
then
$$
|E_j|^{\alpha-\delta}=\left|\bigcup_{i=m}^{m+k}\Delta_{\alpha_1\cdots\alpha_{n_j-1}i}\right|^{\alpha-\delta}=
|\Delta_{\alpha_1\cdots\alpha_{n_j-1}}|^{\alpha-\delta}\left(\sum_{i=m}^{m+k}q_i\right)^{\alpha-\delta}
$$
$$\geq|\Delta_{\alpha_1\cdots\alpha_{n_j-1}}|^{\alpha-\delta}\sum_{i=m}^{m+k}q_i^{\alpha}\geq\sum_{i=m}^{m+k}|
\Delta_{\alpha_1\cdots\alpha_{n_j-1}i}|^{\alpha}.$$
Therefore,
there exists a constant $C$ such that for any $\alpha\in(0,1)$,
$\delta\in(0,\alpha)$, $E\subset[0,1]$,
the following inequality holds:
$$H^{\alpha}(E,\mathcal{A})\leq H^{\alpha}(E,\Phi)\leq CH^{\alpha-\delta}(E,\mathcal{A}).$$
Hence, ${\dim}_{H}(E,\Phi)\leq{\dim}_{H}(E,\mathcal{A})+\delta,\forall\delta\in(0,\alpha)$,
which proves that $${\dim}_H(E)={\dim}_{H}(E,\mathcal{A})={\dim}_{H}(E,\Phi),\forall E\subset[0,1].$$
\smallskip

\textbf{Necessity.}
Suppose inequality (\ref{equ1}) does not hold, i.e.,  we can find $\alpha\in(0,1)$,
$\delta\in(0,\alpha)$ such that for any $N\in\mathbb{N}$ there exists $n,M>N$,
\begin{equation}\label{equ2}
 \left(\sum_{i=n}^{n+M}q_i\right)^{\alpha-\delta}<\sum_{i=n}^{n+M}q_{i}^{\alpha}.
\end{equation}

Now we will construct a cantor-like set $K$ such that $${\dim}_H(K)<{\dim}_H(K,\Phi)$$
(actually, we show that ${\dim}_H(K)\leq\frac{\delta}{2}$ and ${\dim}_H(K,\Phi)\geq\delta$).
\medskip

Let $L<1$ be a positive real number, and $\{\epsilon_{k}\}_{k\geq 1}$ be a decreasing sequence of positive reals with limit zero.

Step 1.
We take $k_1>N$ sufficiently large such that
$$
\left(\sum_{i=k_1}^{\infty}q_i\right)^{\frac{\delta}{2}}\leq L \ \ \text{ and }\ \ \sum_{i=k_1}^{\infty}q_i\leq \epsilon_1.
$$

By inequality (\ref{equ2}),
there exists an integer $M_1>N$ with
$$
\left(\sum_{i=k_1}^{k_1+M_1}q_i\right)^{\alpha-\delta}<\sum_{i=k_1}^{k_1+M_1}q_{i}^{\alpha}.
 $$
Then we choose the following finite union of cylinders of rank $1$,
$$\mathcal{F}_1:=\left\{\bigcup_{\alpha_1=k_1}^{ k_1+M_1}\Delta_{\alpha_1}\right\}.$$

Step $n\geq 1$. Suppose $k_1,k_2,\cdots,k_n$ and $M_1,M_2,\cdots,M_n$ have been defined.
We let $k_{n+1}$ be sufficiently large such that
$$
\sum_{i=1}^{n}\sum_{\alpha_i=k_i}^{k_i+M_i}
\left|\bigcup_{\alpha_{n+1}=k_{n+1}}^{\infty}\Delta_{\alpha_1\cdots\alpha_n\alpha_{n+1}}\right|^{\frac{\delta}{2}}\leq L
\ \ \text{ and } \sum_{i=k_{n+1}}^{\infty}q_i\leq\epsilon_{n+1}.
$$
Then there exists an $M_{n+1}>N$ with
$$
\left(\sum_{i=k_{n+1}}^{k_{n+1}+M_{n+1}}q_i\right)^{\alpha-\delta}<\sum_{i=k_{n+1}}^{k_{n+1}+M_{n+1}}q_{i}^{\alpha}.
$$

Do the same operation step by step,
we get sequences $\{k_n\},\{M_n\}$ of integers.

Define
$$
\mathcal{F}_n:=\left\{\bigcup_{i=k_n}^{k_n+M_n}\Delta_{\alpha_1\cdots\alpha_{n-1}i}\colon
k_j\leq\alpha_{j}\leq k_j+M_j, 1\leq j\leq n-1 \right\}.
$$

Let
$$
K:=\bigcap_{n=1}^{\infty}\bigcup_{\Delta\in\mathcal{F}_n}\Delta.
$$

To this end,
we show that $${\dim}_H(K)\leq\frac{\delta}{2} \ \ \ \text{and}\ \ \  {\dim}_{H}(K,\Phi)\geq\delta.$$

From the construction of $K$,
$\mathcal{F}_n$ covers $K$ for each $n$,
and the $\frac{\delta}{2}$-volume of the coverings $\mathcal{F}_n$ does not exceed $L$,
i.e.,
$$
V_n(\frac{\delta}{2})=\sum_{\Delta\in\mathcal{F}_n}|\Delta|^{\frac{\delta}{2}}
=\sum_{i=1}^{n-1}\sum_{\alpha_i=k_i}^{k_i+M_i}\left|\bigcup_{\alpha_n=k_n}^{k_n+M_n}
\Delta_{\alpha_1\cdots\alpha_{n-1}\alpha_{n}}\right|^{\frac{\delta}{2}}\leq L.
$$
Therefore, $H^{\frac{\delta}{2}}(K)\leq L$, thus ${\dim}_H(K)\leq\frac{\delta}{2}$.

In order to prove ${\dim}_{H}(K,\Phi)\geq\delta$, we will define a mass distribution on $K$.

Let $\xi_n, n\geq 1$ be independent random variables taking values $$k_n, k_n+1,\cdots,k_n+M_n$$
with probabilities
$$\frac{1}{\gamma_n}q_{k_n}^{\alpha},\frac{1}{\gamma_n}q_{k_n+1}^{\alpha},
\cdots,\frac{1}{\gamma_n}q_{k_n+M_n}^{\alpha}$$
respectively,
where $\gamma_n$ is the normalised constants determined by
$$\sum_{i=k_n}^{k_n+M_n}\frac{1}{\gamma_n}q_i^{\alpha}=1.$$
Then by inequality (\ref{equ2}),
we have $$\gamma_n=\sum_{i=k_n}^{k_n+M_n}q_i^{\alpha}>\left(\sum_{i=k_n}^{k_n+M_n}q_i\right)^{\alpha-\delta}.$$

Let $\xi$ be a random variable with independent $Q_{\infty}$-digits,
i.e., $\xi$ is of the form
 $$\xi=\Delta_{\xi_1\xi_2\cdots\xi_n\cdots}.$$
Let $\mu_{\xi}$ be the probability measure of the corresponding random variable $\xi$.
Then
$$
\mu_{\xi}(\Delta_{\alpha_1\cdots\alpha_{n}})=\prod_{i=1}^{n}\frac{1}{\gamma_i}q_{\alpha_i}^{\alpha},
\text{ for } k_i\leq\alpha_i\leq k_i+M_i,1\leq i\leq n.
$$
Recalling that $$|\Delta_{\alpha_1\cdots\alpha_n}|=\prod_{i=1}^{n}q_{\alpha_i},$$
so for any $x\in K$ and $t\in(0,\delta)$
$$
\frac{\mu_{\xi}(\Delta_{\alpha_1(x)\cdots\alpha_n(x)})}{|\Delta_{\alpha_1(x)\cdots\alpha_n(x)}|^t}
=\prod_{i=1}^n\frac{1}{\gamma_i}q_{\alpha_i}^{\alpha-t}
\leq
\prod_{i=1}^n\frac{1}{(\sum_{j=k_i}^{k_i+M_i}q_j)^{\alpha-\delta}}q_{\alpha_i}^{\alpha-t}
\leq\prod_{i=1}^{n}q_{\alpha_i}^{\delta-t}\leq\epsilon_n^{\delta-t}.
$$
So
$$\lim\frac{\mu_{\xi}(\Delta_{\alpha_1(x)\cdots\alpha_n(x)})}{|\Delta_{\alpha_1(x)\cdots\alpha_n(x)}|^t}=0$$
for any $t\in(0,\delta)$ and $x\in K$,
therefore ${\dim}_H(K,\Phi)\geq\delta$ which proves the theorem.


\end{document}